\documentclass[journal]{IEEEtran}
\IEEEoverridecommandlockouts
\usepackage{cite}
\usepackage{amsmath,amssymb,amsfonts}
\usepackage{graphicx}
\usepackage{textcomp}
\usepackage{pgfplots}
\usetikzlibrary{matrix,positioning}
\usepackage{tikz}
\usepackage{threeparttable}
\usepackage{multirow}
\usepackage{array}
\usepackage{booktabs}

\def\BibTeX{{\rm B\kern-.05em{\sc i\kern-.025em b}\kern-.08em
    T\kern-.1667em\lower.7ex\hbox{E}\kern-.125emX}}

\begin{document}
\markboth{This work has been accepted by the 44th Annual Conference of the IEEE Industrial Electronics Society (IECON 2018). This is a preprint.~~~~}%
{Shell \MakeLowercase{\textit{et al.}}: Bare Demo of IEEEtran.cls for IEEE Journals}
\title{Model Predictive Control of H5 Inverter for Transformerless PV Systems with Maximum Power Point Tracking and Leakage Current Reduction\\
}

\author{
	\IEEEauthorblockN{Abdulrahman J. Babqi\IEEEauthorrefmark{1}, Zhehan Yi\IEEEauthorrefmark{2}, Di Shi\IEEEauthorrefmark{2},
	Xiaoying Zhao\IEEEauthorrefmark{2}
	}\\
    \IEEEauthorblockA{\IEEEauthorrefmark{1}Department of Electrical Engineering, Taif University, Taif, Saudi Arabia} \\
	\IEEEauthorblockA{\IEEEauthorrefmark{2}GEIRI North America, San Jose, CA 95134}\\

	Emails: ajbabqi@tu.edu.sa, \{zhehan.yi, di.shi, xiaoying.zhao\}@geirina.net
	\thanks{This work is supported by the SGCC Science and Technology Program Distributed High-Speed Frequency Control under UHVDC Bipolar Blocking Fault Scenario.}		
}

\maketitle

\begin{abstract}
Transformerless grid-connected solar photovoltaic (PV) systems have given rise to more research and commercial interests due to their multiple merits, e.g., low leakage current and small size. In this paper, a model-predictive-control (MPC)-based strategy for controlling transformerless H5 inverter for single-phase PV distributed generation system is proposed. The method further reduces the PV leakage current in a cost-effective and safe manner and it shows a satisfactory fault-ride-through capability. Moreover, for the first of its kind, PV maximum power point tracking is implemented in the single-stage H5 inverter using MPC-based controllers. Various case studies are carried out, which provide the result comparisons between the proposed and conventional control methods and verify the promising performance of the proposed method. 
\end{abstract}


\section{Introduction}
\begin{figure}
	\centering
	\includegraphics[width=1\columnwidth]{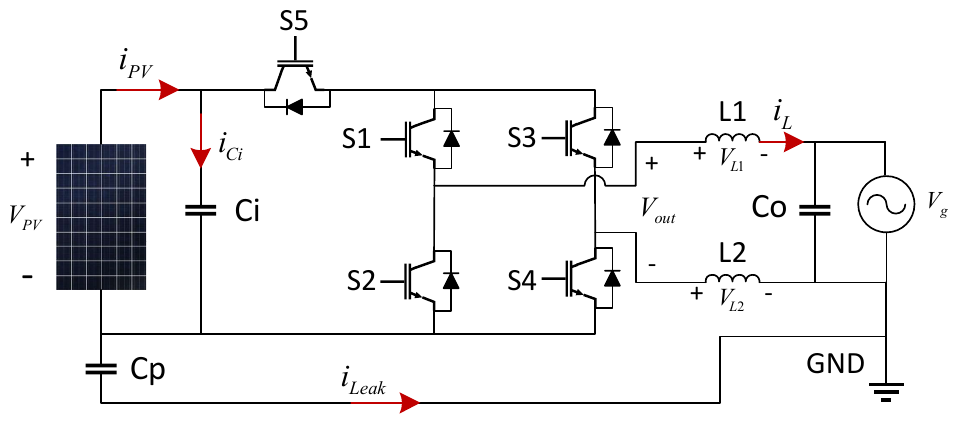}
	\caption{A typical configuration of the H5 transformerless inverter in a PV system.}\label{h5cirucit}
\end{figure}

Renewable distributed energy resources (DERs), such as solar photovoltaic (PV) and wind power systems, have been getting more attentions recently to be used as alternatives of fossil fuels \cite{di, yishen, ahmed}. PV power systems are considered as one of the most attractive renewable DER technologies thanks to the abundance of solar energy and the declining capital and operational expenses \cite{TIE,TSG1}. Generally, PV systems can be interfaced with the utility grid through transformer-isolation or transformerless configurations. Since line frequency transformers are heavy, inefficient, and cost-ineffective for PV systems, transformerless configurations are attracting more and more interests from both research and commercial points of view. However, the lack of galvanic isolation in the transformerless configurations will lead to a common-mode (CM) leakage current between the PV panels and the ground through parasitic capacitors, which reduces the overall efficiency and grid current quality and may cause serious electromagnetic interference and insecurity issues\cite{3}. The parasitic capacitance is approximately 60\,nF to 110\,nF every kilowatt of the PV array \cite{SMA}. Therefore, various inverter topologies with specific modulation strategies have been introduced to suppress the leakage current, in which only a few topologies have been developed into industrial products, e.g., H5, H6, and HERIC inverters. The H5 structure is adopted by the SMA Solar Technology due to its simple topology with the least number of switches \cite{11}. Fig. \ref{h5cirucit} illustrates a typical grid-connected transformerless PV system using a H5 inverter, where $i_{Leak}$ stands for the CM leakage current, C\textsubscript{P} is the PV parasitic capacitor mentioned previously.


\begin{figure*}
	\centering
	\includegraphics[width=1\textwidth]{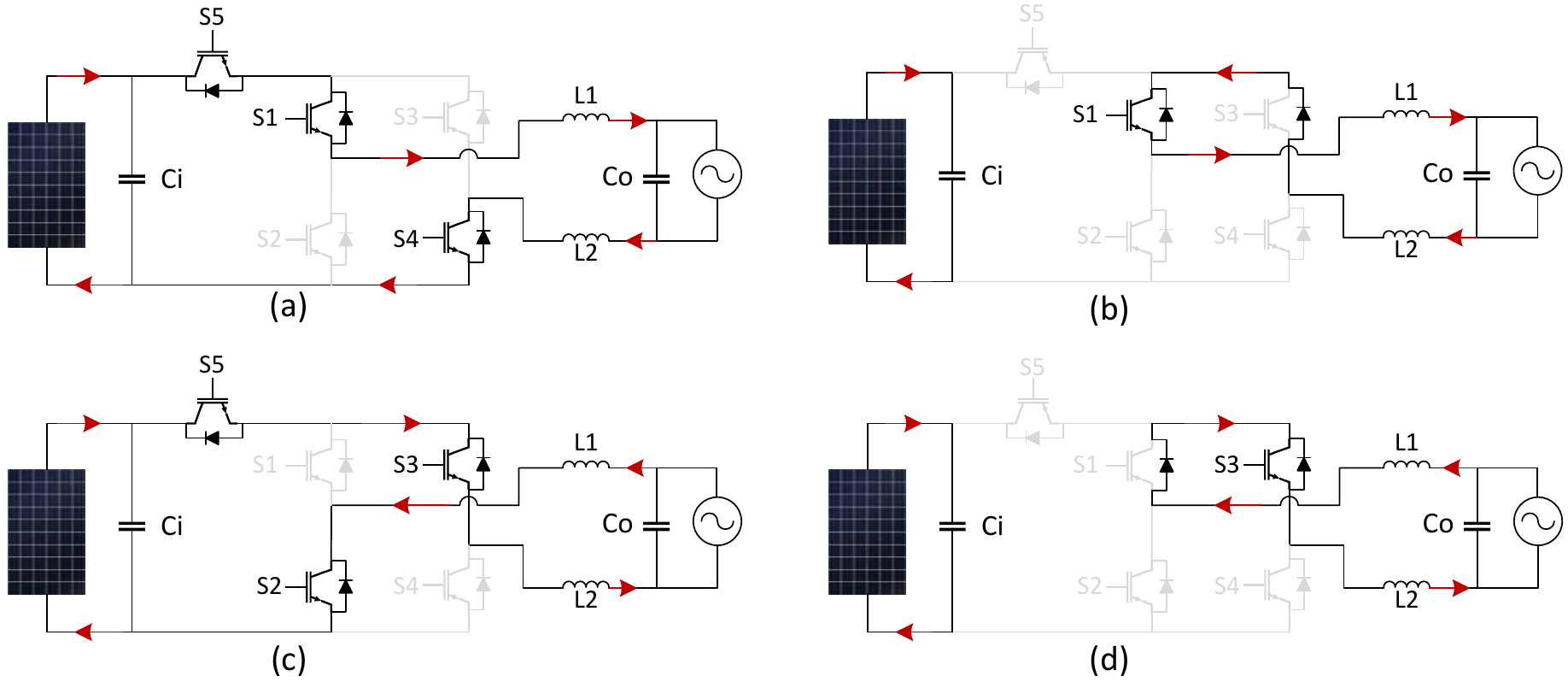}
	\caption{The four operation modes of H5 inverter for grid-connected PV systems.}\label{h5modes}
\end{figure*}  

In order to extract the maximum power of a PV array under different ambient conditions (irradiance and temperature), maximum power point tracking (MPPT) algorithms, such as Perturb \& Observe (P\&O) and Incremental Conductance (In-Cond), are employed to control the power-electronics stage. Although there are numerous existing methods to implement MPPT for grid-tied PV systems, most of them use two-stage cascaded DC/DC-DC/AC converting systems \cite{TSG2, cascaded1,cascaded3} or single-stage DC/AC inverters \cite{inverter1,inverter2,inverter3, inverter4} with PI-based controllers or their variants. These methods may suffer from one or multiple of the following major drawbacks:

\begin{itemize}
	\item PI-based controllers require iterative tuning efforts when system parameters change;
	\item It is relatively difficult to find optimal gain and time constants for the controllers; 
	\item Extra pulse width modulation (PWM) modules are required; 
	\item Some of the methods require multiple stages of costly converters, which reduces the converting efficiency; and
	\item CM leakage current is not considered in most methods.
\end{itemize} 
   

\noindent Furthermore, very few research has tried to control the PV MPPT using a transformerless single-stage H5 inverter.

Model predictive control (MPC) is an optimal control approach which uses the system model and measurements to predict the future behavior of the controlled states based on minimizing a cost function \cite{9, babqi1, babqi2}. It is a fast, robust, and accurate controller that requires little tuning efforts. Recently, MPC has been seen in the literature for controlling PV power systems \cite{6,7,8,yipes2018}. However, these works only design MPC to control DC/DC converters. There is no existing work aiming at designing MPC for H5 inverters. To fill this gap and as an attempt to address the aforementioned issues, this paper proposes a MPC-based strategy for controlling the single-stage transformerless H5 inverter for PV distributed generation systems. The control strategy further reduces the PV leakage current comparing with conventional control methods. Moreover, fast and accurate MPPT is implemented for the single-stage H5 transformerless inverter using MPC. Additionally, the proposed method improves the robustness and the fault-ride-through capability of transformerless PV systems. The rest of the paper is organized as follows: Section II presents the design of H5 inverter and its operation modes; the proposed MPC for H5 inverters is descried in Section III; case studies are carried out in Section IV to verify the control scheme; Section IV concludes the paper.

\section{H5 Inverter Operation Modes}

\begin{figure}
	\centering 
	\includegraphics[width=0.9\columnwidth]{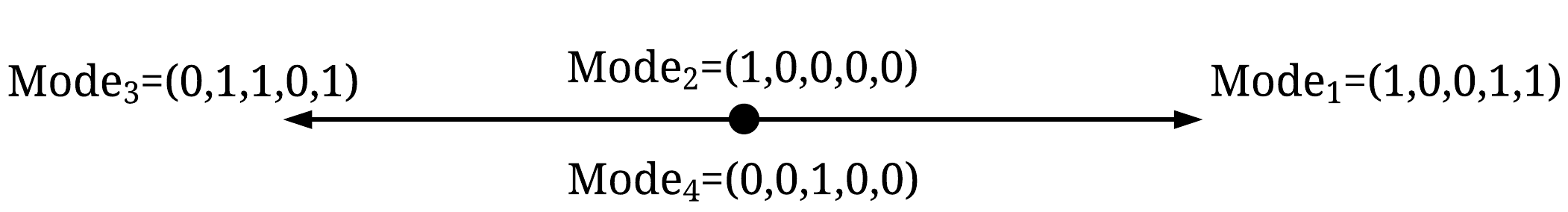}
	\caption{H5 inverter space vector modulation (SVM).}
	\label{H5SVM}
\end{figure}

The topology of a H5 inverter is similar to the single-phase full-bridge inverter by adding an extra DC-bypass switch ``S5'' that disconnects the PV array from the utility grid during the current-freewheeling periods. Fig. \ref{h5cirucit} shows the topology of H5 inverter with the leakage current ($i_{Leak}$) between the PV array and ground. In general, there are four operation modes for H5 inverters, which are depicted Fig. \ref{h5modes}. The first operation mode (Fig. \ref{h5modes}(a)) is the active mode which occurs during the positive-half cycle, where the switches S1, S4, and S5 are conducting and the current flows through S1 and S5 and then returns to the cathode of the PV array through S4. The second mode of operation shown in Fig. \ref{h5modes}(b) is also referred to as the current-freewheeling mode with the zero voltage vector. In this mode, S1 is triggered on, while S4 and S5 are turned-off. The current is conducting through the freewheeling diode of S3. Fig. \ref{h5modes}(c) illustrates the third mode of operation of H5 inverter, which is the active mode that occurs during the negative-half cycle. During mode 3, switches S2, S3, and S5 conduct and the current flows through the inductors L1 and L2 in the opposite direction of that in mode 1. The fourth mode is the freewheeling mode during the zero voltage vector where S2 and S5 are turned-off and S3 is on. Similar to S3 in mode 2 (Fig. \ref{h5modes}(b)), S1 works as a freewheeling diode in mode 4. Table \ref{table1} and Fig. \ref{H5SVM} show the operation modes and space vector modulation (SVM) of the H5 inverter. 
\begin{table}
	\centering
	\caption{H5 inverter Switching states}
	\label{table1}
	\smallskip
	\begin{tabular}{ccccccc}
		\toprule 
		\textbf{Mode} & \textbf{S1} & \textbf{S2} & \textbf{S3} & \textbf{S4} & \textbf{S5} & \textbf{${V_{out}}$}     \\ 
		\midrule
		1 & 1 & 0 & 0 & 1 & 1 & $V_{PV}$ \\ 
		2 & 1 & 0 & 0 & 0 & 0 & 0 \\ 
		3 & 0 & 1 & 1 & 0 & 1 & $-V_{PV}$ \\
		4 & 0 & 0 & 1 & 0 & 0 & 0 \\ 
		\bottomrule
	\end{tabular}
\end{table}

\section{The Proposed Model Predictive Control for Transformerless H5 Inverters Enabling PV MPPT}

\subsection{State Predictions of H5 inverter}
As is mentioned above, in the first mode of operation, S1, S4, and S5 are conducting. The system model can be derived by KCL and KVL, respectively:  
\begin{equation} 
i_{PV} = i_{Ci} + i_{L}
 \label{mode1}
\end{equation}
\begin{equation} 
V_{PV} = V_{L1} + V_{L2} + V_{g}
\end{equation} 

\noindent where $i_{PV}$ and $V_{PV}$ are the PV array output current and voltage, respectively. $i_{Ci}$ and $i_L$ are the currents through capacitor Ci and the inductor L1, respectively. $V_{L1}$ and $V_{L2}$ are the voltages across inductor L1 and L2 and $V_g$ is the utility grid voltage. 

In the second mode of operation where only S1 is turned-on, the system model is:

\begin{equation} 
i_{PV} = i_{Ci},
\label{mode2} 
\end{equation} 
\begin{equation} 
V_{PV} = V_{Ci}. 
\end{equation} 

\noindent where $V_{Ci}$ is the capacitor Ci voltage. 

In the third mode, S2, S3, and S5 are closed and the system model can be written as:

\begin{equation} 
i_{PV} = i_{Ci} + i_{L},
 \label{mode3}
\end{equation}
\begin{equation} 
V_{PV} = V_{L1} + V_{L2} + V_{g}.
\end{equation} 

\noindent During the fourth mode, the only closed switch is S3 and the system model can be given as:
\begin{equation} 
i_{PV} = i_{Ci},
 \label{mode4}
\end{equation}
\begin{equation} 
V_{PV} = V_{Ci}. 
\end{equation} 

\noindent Additionally, the capacitor current and inductor voltages can be expressed as:
\begin{equation} 
i_{Ci} = \frac{dV_{PV}}{dt},\\
V_{L1} = \frac{di_{L1}}{dt},\\
V_{L2} = \frac{di_{L2}}{dt}.\\
\end{equation}

For a PV system, the MPPT is realized by forcing the PV operating pointing to be around maximum power point, namely, to control $V_{PV}$ to track the MPPT reference voltage $V_{ref}$. Therefore, we have to predict the future values of PV array voltage, $V_{PV}(k+1)$ per horizon step. To this end, (\ref{mode1}), (\ref{mode2}), (\ref{mode3}), and (\ref{mode4}) must be discretized. Using the forward finite difference formula for the derivative 
\begin{equation}
\frac{dx}{dt}\approx\frac{x(k+1)-x(k)}{T_s}, 
\end{equation}

\noindent where $T_s$ is the sampling period, the future values of the PV output voltage for the aforementioned four operation modes are predicted as:
\begin{equation} 
V_{PV1}(k+1) = V_{PV}(k) + \frac{T_s}{Ci} [i_{PV}(k) - i_{L}(k)]
\end{equation} 
\begin{equation} 
V_{PV2}(k+1) = V_{PV}(k) + \frac{T_s}{Ci}~i_{PV}(k)
\end{equation} 
\begin{equation} 
V_{PV3}(k+1) = V_{PV}(k) + \frac{T_s}{Ci} [i_{PV}(k) - i_{L}(k)]
\end{equation} 
\begin{equation} 
V_{PV4}(k+1) = V_{PV}(k) + \frac{T_s}{Ci}~i_{PV}(k)
\end{equation}

\noindent It is noteworthy that for the PV array voltage, the predictions of mode 1 and 3 are identical, so are mode 2 and 4. This will reduce the mode switching frequency and thus the CM leakage current, which will be seen in the verification in the case studies (Section IV).

\begin{figure}
	\centering
	\includegraphics[width=1\columnwidth]{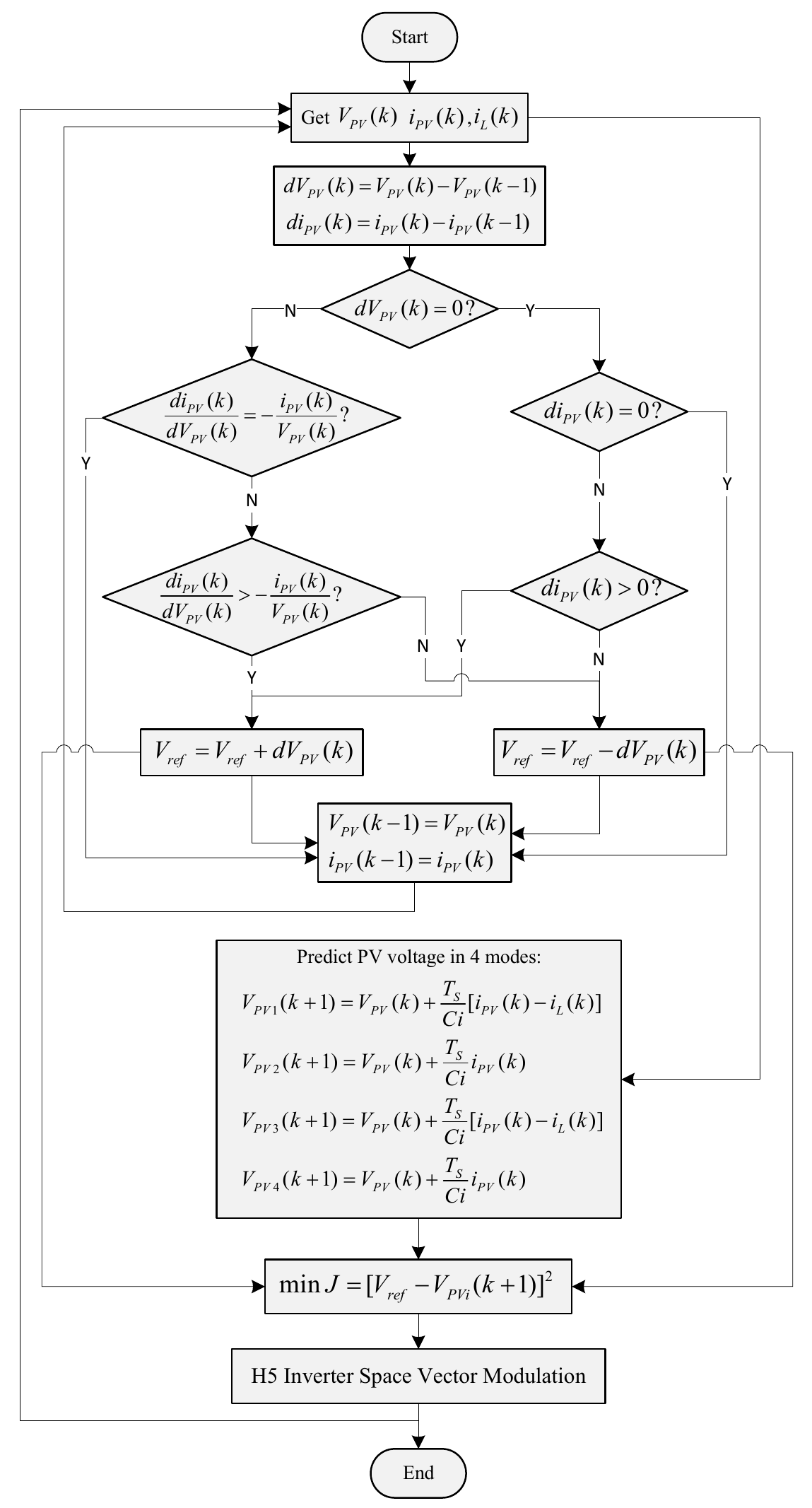}
	\caption{The proposed MPC-based PV MPPT control for transformerless H5 inverters.}\label{flowchart}
\end{figure}

\subsection{Maximum Point Point Control of H5 Inverter Using MPC}
The proposed method is modified from one of the most popular MPPT algorithms, the Incremental Conductance method \cite{yi_phd}, for extracting the maximum power output of the transformerless PV array with H5 inverter under varying irradiance and temperature. Fig. \ref{flowchart} depicts the process of the proposed algorithm in detail. At every sampling period Ts, the controller samples the values of $V_{PV}$, $i_{PV}$, and $i_L$ from the transformerless PV system and follows the procedures of the flowchart to determine the optimal control inputs. After predicting the future value of the PV output voltage, a quadratic cost function:

\begin{equation} 
J = [V_{ref} - V_{PV}(K+1)]^2
\label{costfunction}
\end{equation}  

\noindent is designed to quantify the difference between the MPPT reference voltage $V_{ref}$ and the future PV array voltage $V_{PV}(k+1)$. The control process is then transferred into a optimization problem, which minimizes the quadratic function (\ref{costfunction}), i.e., the error between the reference and predicted value, and select the optimal control inputs with the least cost. This is achieved by evaluating each possible scenario (i.e., the four modes of operation) and selecting the best operation mode at every sampling step. Once the optimal operation mode is determined, appropriate gating signals are sent to the H5 inverter switches. 

Via optimization, the controller will automatically select the switch signals that lead to a minimum error between the controlling states and references, which eliminates the tuning efforts that required by conventional controllers. Moreover, the switching signals will be directly applied to the H5 inverter without the needs for an extra PWM module, which lowers the cost and complexity of the control system.

\section{Case Studies}
To examine the performance of the proposed control strategy for H5 inverters, multiple case studies are carried out in this section. The transformerless PV system with the same configuration in Fig. \ref{h5cirucit} is modeled in the PSCAD/EMTDC platform, while the proposed algorithm is implemented using Fortran. The numerical values of the tested system parameters are provided in Table \ref{table2}. The parameters of the PV array are measured under standard testing condition (STC, irradiance=1000\,W/m\textsuperscript{2}, temperature=25\textdegree C) It is noteworthy that, although the proposed method is verified in a testbed with certain parameters, it is scalable to work for different transformerless PV systems. For systems with other configurations, similar approach can be used to design the controller. The tested cases are elaborated below.

\subsection*{Case Study 1}

\begin{figure}
	\centering 
	\includegraphics[width=1\columnwidth]{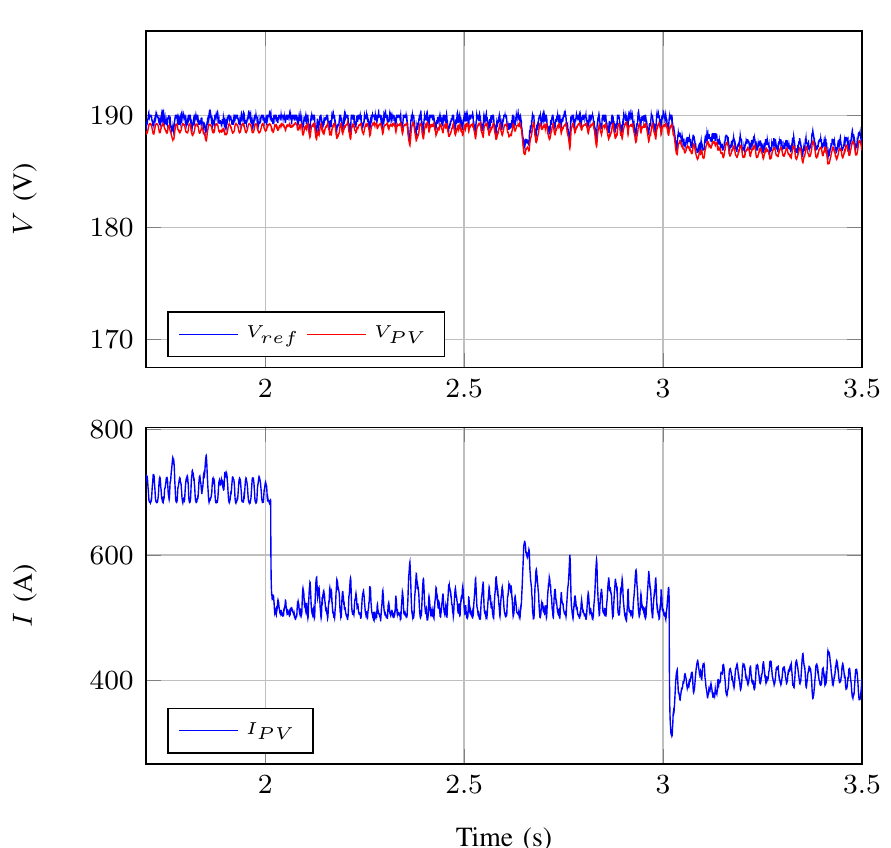}
	\caption{PV output voltage and current using the proposed control method.}
	\label{fig.5} 
\end{figure}

This case verifies the MPPT capability of the proposed control method, which aims at extracting the maximum power of the PV array under varying irradiance situations. Fig. \ref{fig.5} shows the PV array instantaneous output voltage ($V_{PV}$) and current ($I_{PV}$) as well as the MPPT voltage reference ($V_{ref}$). At the beginning, the irradiance is set to 1000 W/m\textsuperscript{2} and the temperature is $25^\circ$C. At $\textit{t}=2~s$, the irradiance drops from 1000 W/m\textsuperscript{2} to 800 W/m\textsuperscript{2}. It can be seen from Fig. \ref{fig.5} that the output current of the PV array decreases as the irradiance changes without any undershooting. Moreover, the PV output voltage is tracking its reference $V_{ref}$ closely. A further decline of the irradiance occurred at $\textit{t}=3~s$. The irradiance reduces from 800 W/m\textsuperscript{2} to 600 W/m\textsuperscript{2}. It is demonstrated that both PV voltage and current change their values correspondingly because of the new irradiance level. Again, it is clear that the PV output voltage follows its reference value. Therefore, the proposed control strategy provides a fast response as well as good dynamic performance under the varying irradiance conditions.

\begin{table}
	\centering
	\caption{Case Study System Parameters}
	\label{table2}
	\smallskip
	\begin{tabular}{lcc}
		\toprule 
		\textbf{Parameter} & {\textbf{Symbol}} & {\textbf{Value}} \\ 
		\midrule
		Standard Testing Irradiance  & G          & 1000\,W/m\textsuperscript{2}\\
		Standard Testing Temperature & T          & 25\textdegree C \\ 
		PV Array Maximum Power (STC)   & $P_{max}$  & 133\,kW \\ 
		PV Array Maximum Point Point Voltage (STC)  & $V_{ref}$  & 190\,V \\ 
		DC Capacitor  & Ci & 5000\,uF   \\
		Output Filter Inductor 1  &  L1 &  1\,mH \\
		Output Filter Inductor 2  &  L2 &  1\,mH \\
		Output Filter Capacitor   &  Co &  5000\,uF \\
		Grid Side Voltage & $V_g$ &  100\,V  \\
		Sampling Period & T\textsubscript{s} &  100\,us  \\
		PV Parasitic Capacitance & C\textsubscript{P} &  13300\,nF  \\
		\bottomrule
	\end{tabular}
\end{table}

\subsection*{Case Study 2} 

\begin{figure}
	\centering 
	\includegraphics[width=1\columnwidth]{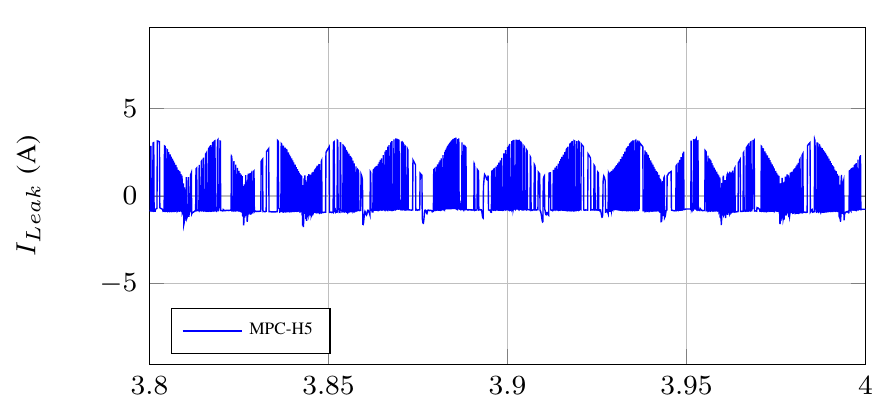}
	\includegraphics[width=1\columnwidth]{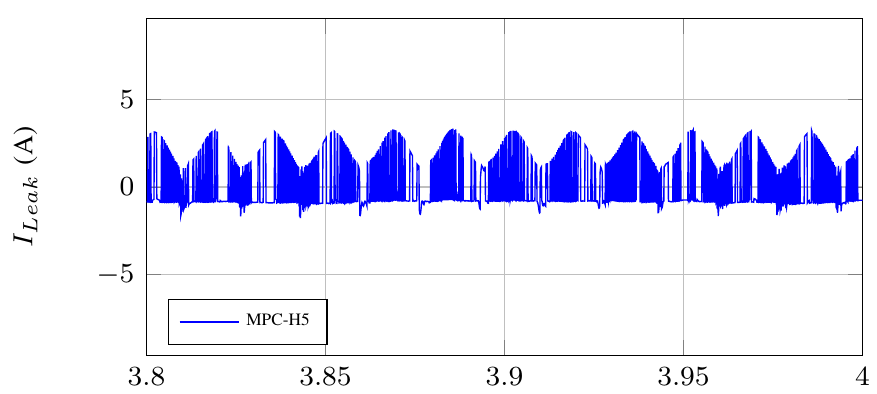}
	\includegraphics[width=1\columnwidth]{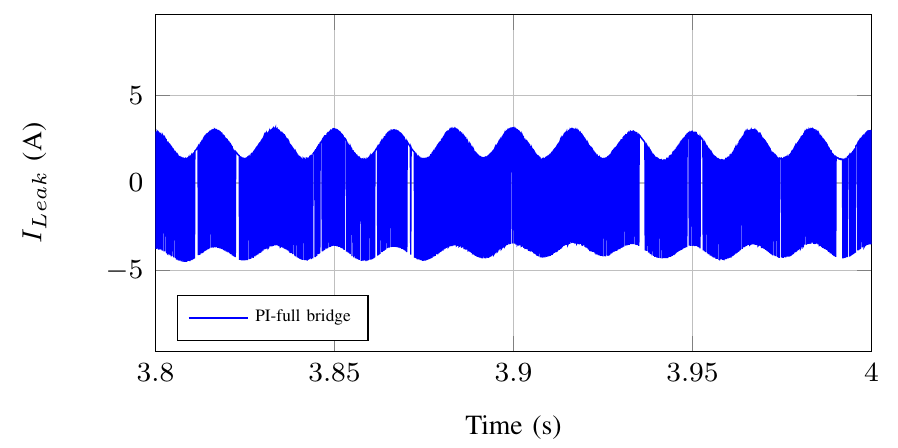}
	\caption{Leakage current comparison among the MPC-controlled H5 inverter, PI-controlled H5 inverter, and PI-controlled full-bridge inverter.}
	\label{fig.6} 
\end{figure}

\begin{figure}
	\centering 
	\includegraphics[width=1\columnwidth]{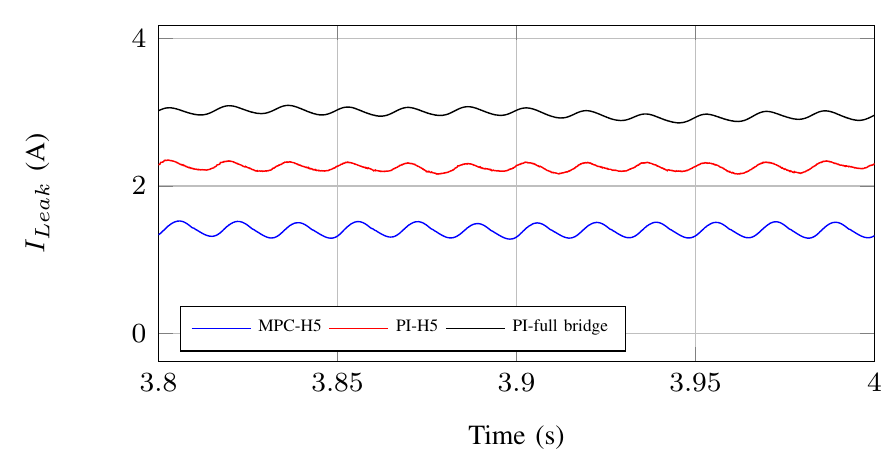}
	\caption{RMS values of the leakage current}
	\label{fig.7} 
\end{figure}

The following case aims at validating the capability of the proposed control strategy for further reducing the CM leakage current in a transformerless PV system with a H5 inverter. Fig. \ref{fig.6} illustrates the instantaneous CM leakage currents of the transformerless PV system by the proposed MPC and conventional PI controllers, as well as the leakage current of a PV system with a single-phase full-bridge inverter by a PI controller. Fig. \ref{fig.7} presents their RMS values. To present a reasonable comparison, the MPPT algorithm for all these case are based on the Incremental Conductance method. From these figures, it is obvious that H5 inverter itself reduces the CM leakage current (RMS) from 3\,A to 2.5\,A approximately (black and red curves in Fig. \ref{fig.7}). 

\begin{figure}
	\centering 
	\includegraphics[width=1\columnwidth]{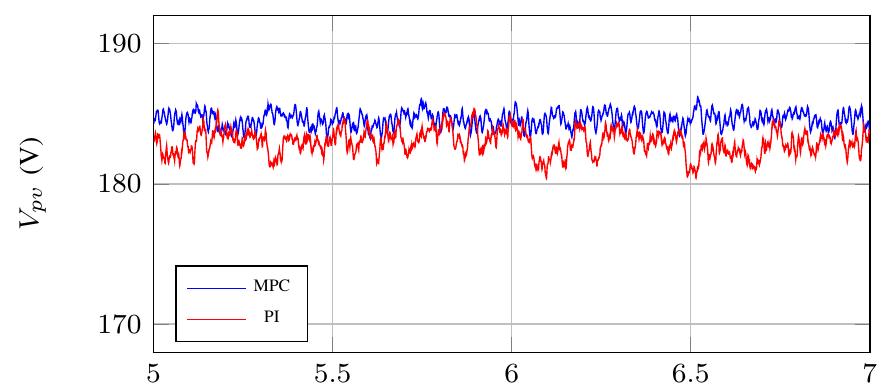}
	\includegraphics[width=1\columnwidth]{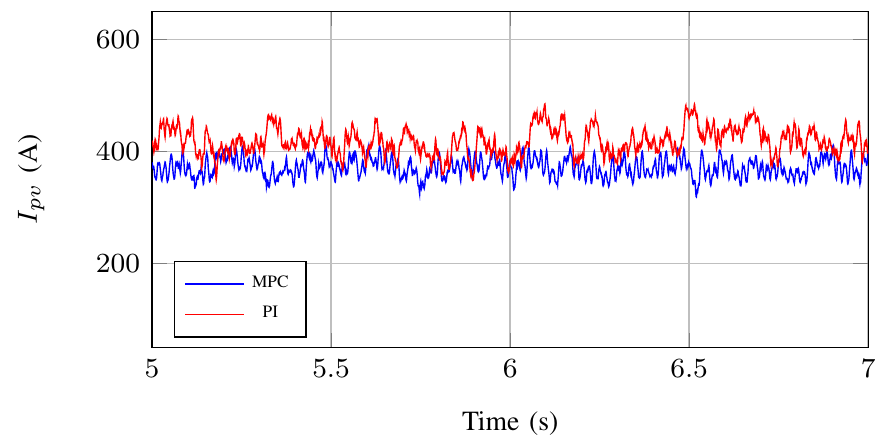}
	\caption{Comparison of the PV output voltage and current using the proposed controller and a conventional PI controller.}
	\label{fig.8} 
\end{figure}

The results also presents that, for the same configuration (transformerless H5 inverter), the proposed method further reduces the leakage current by almost 50\% compared with conventional PI controller (blue curve in Fig. \ref{fig.7}). This is because that, during operation, the proposed method reduces the switching modes of the H5 inverter as is analyzed in Section III. The large leakage current that occurs with the conventional controller will affect the reliability and efficiency of the PV system. More importantly, it may cause safety hazards to the system operator and maintenance personnel. Fig. \ref{fig.8} demonstrates the PV output voltage and current using both the proposed control strategy and conventional PI control for H5 inverter, which shows that the leakage current affects the performance of the PV system and makes it more oscillatory while proposed method gives a smoother and steadier performance.

\begin{figure}
	\centering 
	\includegraphics[width=1\columnwidth]{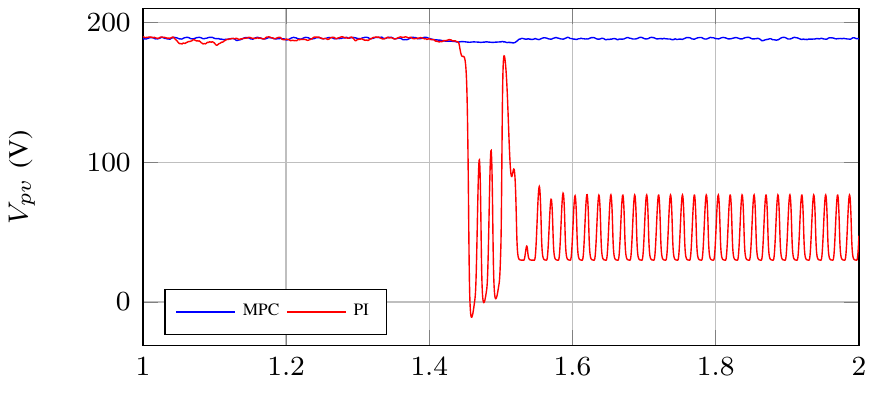}
	\includegraphics[width=1\columnwidth]{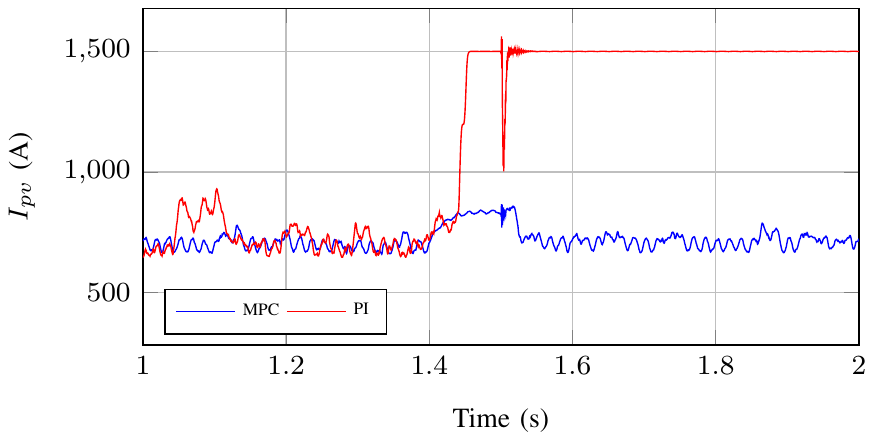}
	\caption{Output PV voltage and current using both proposed controller and conventional PI.}
	\label{case3} 
\end{figure}

\subsection*{Case Study 3} 

This case demonstrates the fault-ride-through capability of the proposed control method. To this end, a ground fault was applied to system at the output terminal of the PV array. Fig. \ref{case3} plots the PV output voltage and current of H5 inverter (using both the proposed method and PI controller) before, during, and after the fault. The ground fault is applied at $\textit{t}=1.4~s$ and it is cleared after 100 ms (Fig. \ref{case3}). It can be seen that the system controlled by the PI controller is vulnerable. The PV voltage and current become unstable during and even after the clearance of fault (red curves). The PV current increases immediately when the fault occurs, which may cause damage to the system if no further protection actions are applied. Nevertheless, the proposed control strategy shows a robust performance (blue curves) and a better fault-ride-through capability under faulted conditions. The PV voltage tracks its reference value $V_{ref}$ both during and after the fault cleared, while the current is limited within a reasonable range.

\section{Conclusions}

This paper introduces an innovative control strategy for transformerless grid-connected PV systems with H5 inverters. A model-predictive-controlled method is designed to extract the PV maximum power under various operational conditions. The control strategy predicts the future behavior of the PV output voltage and generates the optimal control signals for the H5 inverter, which minimizes the error between the reference and the controlled variable. The case studies verifies that the proposed method provides a better dynamics response comparing with conventional methods. Moreover, the control strategy further reduces the CM leakage current in the H5 inverter by almost 50\% compared with conventional PI controller. It also demonstrates the robustness and fault-ride-through capability of the proposed method.

\bibliographystyle{IEEEtran}

\begin{thebibliography}{10}
	\providecommand{\url}[1]{#1}
	\csname url@samestyle\endcsname
	\providecommand{\newblock}{\relax}
	\providecommand{\bibinfo}[2]{#2}
	\providecommand{\BIBentrySTDinterwordspacing}{\spaceskip=0pt\relax}
	\providecommand{\BIBentryALTinterwordstretchfactor}{4}
	\providecommand{\BIBentryALTinterwordspacing}{\spaceskip=\fontdimen2\font plus
		\BIBentryALTinterwordstretchfactor\fontdimen3\font minus
		\fontdimen4\font\relax}
	\providecommand{\BIBforeignlanguage}[2]{{%
			\expandafter\ifx\csname l@#1\endcsname\relax
			\typeout{** WARNING: IEEEtran.bst: No hyphenation pattern has been}%
			\typeout{** loaded for the language `#1'. Using the pattern for}%
			\typeout{** the default language instead.}%
			\else
			\language=\csname l@#1\endcsname
			\fi
			#2}}
	\providecommand{\BIBdecl}{\relax}
	\BIBdecl
	
	\bibitem{di}
	D.~Shi, X.~Chen, Z.~Wang, X.~Zhang, Z.~Yu, X.~Wang, and D.~Bian, ``A
	distributed cooperative control framework for synchronized reconnection of a
	multi-bus microgrid,'' \emph{IEEE Transaction on Smart Grid}, vol.~PP,
	no.~99, pp. 1--1, 2017.
	
	\bibitem{yishen}
	Y.~Wang, Z.~Yi, D.~Shi, Z.~Yu, B.~Huang, and Z.~Wang, ``Optimal distributed
	energy resources sizing for commercial building hybrid microgrids,''
	\emph{arXiv preprint arXiv:1803.00660}, 2018.
	
	\bibitem{ahmed}
	A.~Aldhaheri and A.~H. Etemadi, ``Stabilization and performance preservation of
	dc-dc cascaded systems by diminishing output impedance magnitude,''
	\emph{IEEE Transactions on Industry Applications}, vol.~54, no.~2, pp.
	1481--1489, March 2018.
	
	\bibitem{TIE}
	Z.~Yi and A.~H. Etemadi, ``Line-to-line fault detection for photovoltaic arrays
	based on multiresolution signal decomposition and two-stage support vector
	machine,'' \emph{IEEE Transactions on Industrial Electronics}, vol.~64,
	no.~11, pp. 8546--8556, Nov 2017.
	
	\bibitem{TSG1}
	Z.~Yi and A.~H. Etemadi, ``Fault detection for photovoltaic systems based on multi-resolution
	signal decomposition and fuzzy inference systems,'' \emph{IEEE Transaction on
		Smart Grid}, vol.~8, no.~3, pp. 1274--1283, May 2017.
	
	\bibitem{3}
	H.~Xiao and S.~Xie, ``Leakage current analytical model and application in
	single-phase transformerless photovoltaic grid-connected inverter,''
	\emph{IEEE Transactions on Electromagnetic Compatibility}, vol.~52, no.~4,
	pp. 902--913, Nov 2010.
	
	\bibitem{SMA}
	``Capacitive leakage currents,'' SMA Solar Technology AG, Ableitstrom-TI-en-25,
	Tech. Rep.
	
	\bibitem{11}
	H.~Li, Y.~Zeng, B.~Zhang, Q.~Zheng, R.~Hao, and Z.~Yang, ``An improved h5
	topology with low common-mode current for transformerless pv grid-connected
	inverter,'' \emph{IEEE Transactions on Power Electronics}, pp. 1--1, 2018.
	
	\bibitem{TSG2}
	Z.~Yi, W.~Dong, and A.~H. Etemadi, ``A unified control and power management
	scheme for {PV}-battery-based hybrid microgrids for both grid-connected and
	islanded modes,'' \emph{IEEE Transactions on Smart Grid}, vol.~PP, no.~99,
	pp. 1--1, 2017.
	
	\bibitem{cascaded1}
	S.~Jiang, D.~Cao, Y.~Li, and F.~Z. Peng, ``Grid-connected boost-half-bridge
	photovoltaic microinverter system using repetitive current control and
	maximum power point tracking,'' \emph{IEEE Transactions on Power
		Electronics}, vol.~27, no.~11, pp. 4711--4722, Nov 2012.
	
	\bibitem{cascaded3}
	C.~Jain and B.~Singh, ``A three-phase grid tied spv system with adaptive dc
	link voltage for cpi voltage variations,'' \emph{IEEE Transactions on
		Sustainable Energy}, vol.~7, no.~1, pp. 337--344, Jan 2016.
	
	\bibitem{inverter1}
	R.~Kadri, J.~P. Gaubert, and G.~Champenois, ``An improved maximum power point
	tracking for photovoltaic grid-connected inverter based on voltage-oriented
	control,'' \emph{IEEE Transactions on Industrial Electronics}, vol.~58,
	no.~1, pp. 66--75, Jan 2011.
	
	\bibitem{inverter2}
	Y.~Yang, H.~Wang, F.~Blaabjerg, and T.~Kerekes, ``A hybrid power control
	concept for pv inverters with reduced thermal loading,'' \emph{IEEE
		Transactions on Power Electronics}, vol.~29, no.~12, pp. 6271--6275, Dec
	2014.
	
	\bibitem{inverter3}
	F.~M. de~Oliveira, S.~A.~O. da~Silva, F.~R. Durand, L.~P. Sampaio, V.~D. Bacon,
	and L.~B.~G. Campanhol, ``Grid-tied photovoltaic system based on pso mppt
	technique with active power line conditioning,'' \emph{IET Power
		Electronics}, vol.~9, no.~6, pp. 1180--1191, 2016.
	
	\bibitem{inverter4}
	M.~Fortunato, A.~Giustiniani, G.~Petrone, G.~Spagnuolo, and M.~Vitelli,
	``Maximum power point tracking in a one-cycle-controlled single-stage
	photovoltaic inverter,'' \emph{IEEE Transactions on Industrial Electronics},
	vol.~55, no.~7, pp. 2684--2693, July 2008.
	
	\bibitem{9}
	A.~J. Babqi and A.~H. Etemadi, ``{MPC}-based microgrid control with
	supplementary fault current limitation and smooth transition mechanisms,''
	\emph{IET Generation, Transmission Distribution}, vol.~11, no.~9, pp.
	2164--2172, 2017.
	
	\bibitem{babqi1}
	A.~J. Babqi, Z.~Yi, and A.~H. Etemadi, ``Centralized finite control set model
	predictive control for multiple distributed generator small-scale
	microgrids,'' in \emph{2017 North American Power Symposium (NAPS)}, Sept
	2017, pp. 1--5.
	
	\bibitem{babqi2}
	A.~J. Babqi, ``Finite control set model predictive control for multiple
	distributed generators microgrids,'' \emph{Ph.D. Dissertation, The George
		Washington University}, 2018.
	
	\bibitem{6}
	M.~Metry, M.~B. Shadmand, R.~S. Balog, and H.~Abu-Rub, ``Mppt of photovoltaic
	systems using sensorless current-based model predictive control,'' \emph{IEEE
		Transactions on Industry Applications}, vol.~53, no.~2, pp. 1157--1167, March
	2017.
	
	\bibitem{7}
	M.~Mosa, M.~B. Shadmand, R.~S. Balog, and H.~A. Rub, ``Efficient maximum power
	point tracking using model predictive control for photovoltaic systems under
	dynamic weather condition,'' \emph{IET Renewable Power Generation}, vol.~11,
	no.~11, pp. 1401--1409, 2017.
	
	\bibitem{8}
	M.~B. Shadmand, R.~S. Balog, and H.~Abu-Rub, ``Model predictive control of pv
	sources in a smart dc distribution system: Maximum power point tracking and
	droop control,'' \emph{IEEE Transactions on Energy Conversion}, vol.~29,
	no.~4, pp. 913--921, Dec 2014.
	
	\bibitem{yipes2018}
	Z.~Yi, A.~Babqi, Y.~Wang, D.~Shi, A.~Etemadi, Z.~Wang, and B.~Huang,
	``Finite-control-set model predictive control (fcs-mpc) for islanded hybrid
	microgrids,'' \emph{arXiv preprint arXiv:1802.04435}, 2018.
	
	\bibitem{yi_phd}
	Z.~Yi, ``Solar photovoltaic ({PV}) distributed generation systems - control and
	protection,'' \emph{Ph.D. Dissertation, The George Washington University},
	2017.
	
\end{thebibliography}


\end{document}